# REJOINDER TO "BREAKDOWN AND GROUPS" BY P. L. DAVIES AND U. GATHER

By P. Laurie Davies and Ursula Gather

*University of Duisburg–Essen and Technical University Eindhoven, and University of Dortmund*

We thank all the discussants for their contributions and in particular we wish to thank Hampel. The concept of breakdown point goes back to his Ph.D. thesis [Hampel (1968)] and he was the first to exhibit a high breakdown equivariant regression estimate now known as the least median of squares [Hampel (1975)], a fact which is sometimes forgotten. These two sources are the starting point of the present discussion. In his contribution Hampel gives us insight into the thoughts which led to his definition of breakdown point, intended as it was to complement the infinitesimal behavior of a functional as described by the influence function. Hampel emphasizes that equivariance considerations were not part of his definition and he had in mind correlation statistics "where there is no equivariance at all." He considers correlation in some detail and, as we disagree with him on this very topic, we give a detailed analysis of correlation statistics in our rejoinder. We hope that this will help clarify the issues involved.

**1. On breakdown.** The first signification of the word "breakdown" given in the Oxford Dictionary starts with the following subsignification:

"**1. a.** The act of breaking and falling down: a ruinous downfall, a collapse."

**2. Breakdown to points and variations.** Genton and Lucas and Oja argue for the usefulness of the breakdown concept in situations not covered by the results of our paper. They claim that at least in an intuitive sense breakdown occurs if the value of a functional is driven to the boundary or to an interior point which is independent of the uncontaminated sample. A formal definition of breakdown point is given which is intended to cover such possibilities. A first version is to be found in Genton and Lucas (2003) and is referred to by Oja. It defines the breakdown point as the smallest amount of contamination which can cause the statistic to assume only a finite number of values independently of the uncontaminated observations. On the basis of this definition the arithmetic mean is claimed to have a finite-sample breakdown point of $1/n$. The argument is as follows: if the first observation of the sample is contaminated, $(\xi_1, y_2, \ldots, y_n)$, and we let $\xi_1$ tend to infinity, then the sample mean tends to infinity, that is, to a single value which is independent of $y_2, \ldots, y_n$. However, for any finite value of $\xi_1$ the arithmetic mean takes on a continuum of values on varying the uncontaminated part of the sample. The only way of reducing the arithmetic mean to a single value is to introduce the symbol $\infty$ as a possible value for the contamination. The symbol $\infty$ is thus elevated to a real entity for data. The new definition avoids this but our reaction is similar: any definition of breakdown based on the concept of Lebesgue measure zero must be at fault. According to the new definition the functional

$$T(P_n) = \max\{-n, \min\{n, T_{\mathrm{LS}}(P_n)\}\}$$

has a breakdown point of $1/n$. We perturb it and put

$$T^*(P_n) = T(P_n) + \frac{1}{n}\int \sin(x)\,dP_n(x).$$

The set of values taken on by $T^*(P_n)$ as we vary the uncontaminated part of the sample has Lebesgue measure at least $2/n^2$ as long as not all the sample is contaminated, and the breakdown point is therefore 1. As the perturbation tends with $1/n$ to zero, $T^*$ remains consistent and asymptotically normal at the model. Oja mentions the classical skewness statistic

$$b_1 = \frac{((1/n)\sum(x_i - \bar{x})^3)^2}{((1/n)\sum(x_i - \bar{x})^2)^3},$$

but this can be treated in the same manner by putting

$$b_1^* = b_1 + \sin(nb_1),$$

which is still invariant but does not converge.

A second criticism we made of the definition of Genton and Lucas (2003) is that any realizable functional immediately breaks down for the simple reason that it can only take on a finite number of values; all data and statistics are of finite precision. Genton and Lucas mention this in their contribution as a weakness of the new definition and so it is. No reasonable definition of breakdown can rely on the myth of a continuum of possible values for a statistic or the associated myth of infinite precision. When applying mathematics to applied problems it is important that the discrete problem can be well approximated by the continuous one. Genton and Lucas' use of infinite precision and a continuum of values and sets of Lebesgue measure zero is not of this sort. Their continuous formulations do not approximate the discrete world of statistics.

We point out further that the definitions of Genton and Lucas, and also of Oja (Definition 4), represent a complete break with the meaning of breakdown as it is used in statistics. Transferred to the statistical context the "ruinous downfall" of Section 1 is expressed in terms of distances and arbitrarily large bias. None of this is present in a concept of breakdown in terms of the number or Lebesgue measure of the set of all possible limits of contaminated samples. No mention is made of bias, that is, how far the value of the statistic can move from its value at the uncontaminated sample for a given amount of contamination. Yet it is this which has motivated robust statistics from the influence function via bias to breakdown point. In a sense the proposal put forward by Genton and Lucas is the very opposite of this. Rather than moving arbitrarily far, the statistic has broken down if it does not move at all. It is said that it then cannot convey any information in the sample. Even this is not always the case. Consider the statistical functional $T_{75}$ which takes on the fixed value of 75 for all data sets. This has a breakdown point of 0 according to the Genton–Lucas definition. German insurance companies are required to use life expectancies specified by law. In the case of a male they could, for example, be forced to use the functional $T_{75}$ to estimate life expectancy in years. The effect can be felt but it is not a ruinous downfall. It would be a ruinous downfall for the German insurance companies if they had to use a value of 65 and the reason is that 65 differs from the experienced lengths of life much more than does the fixed value of 75. Here as in the usual definition of breakdown it is the discrepancy which is important.

**3. Perturbations.** The criticism we gave of Genton and Lucas's definitions of breakdown has wider implications. We regard robust statistics as a perturbation theory for statistics. In particular, robust statistics must concern itself with perturbations of models and data sets and, in consequence, it must be able to deal with finite precision. The perturbations involved should be realistic ones and this will in general exclude perturbations described by the gross error neighborhood, which is simply too small. Unfortunately, the idea of stability under perturbations is sometimes lost, especially in theoretical work. Suppose a theorem on the existence and uniqueness of a functional requires assumptions about the existence and differentiability of a density function. These assumptions should then not be referred to as "under weak assumptions" but rather as "under very restrictive assumptions which violate the spirit of robustness." Densities disappear under perturbations, likelihood disappears under perturbations as does the property of being a Lebesgue set of measure zero, efficiency is pathologically discontinuous, and so on. Perturbations and their consequences should be taken seriously by all who work in the area of robust statistics.

**4. Affine equivariant location functionals.** The example of location functionals makes use of only the translation group although it seems natural to require affine equivariance. The problem is that for the affine group we have $G_1 = \varnothing$ since if we iterate $\mathcal{A}(\theta) = A(\theta) + b$ this will in general not tend to infinity so Theorem 3.1 is not applicable. The highest breakdown point for translation equivariant functionals is $1/2$ but there are affine equivariant location functionals which are based on scatter functionals and which have a breakdown point of at least that of the scatter functional, namely $(1 - \Delta(P))/2$. The gap has not been closed but Rousseeuw gives a sufficient condition for the bound $(1 - \Delta(P))/2$ to hold. His argument makes use of the convex hull which can be seen as a form of scatter functional albeit with a low breakdown point. In Davies and Gather (2002) we showed that the bound $1/2$ is attainable at least at some empirical measures so that the gap remains.

**5. Metrics on $\mathcal{P}$.** We agree with Hampel's comments on the gross error neighborhood but we do not like either of the alternatives he suggests. First, total variation is not much better than the gross error neighborhood; a distribution $Q$ lies in the $\varepsilon$ total variation neighborhood of $P$ if and only if $Q - P = \varepsilon(H_1 - P) - \varepsilon(H_2 - P)$ for some distributions $H_1$ and $H_2$ [see Rieder (2000), page 7]. Second, the Prohorov metric defined by

$$d_{\mathrm{pr}}(P, Q) = \inf\{\varepsilon > 0 : P(A) \leq Q(A^\varepsilon) + \varepsilon\}, \tag{5.1}$$

where

$$A^\varepsilon = \{x : d(x, A) < \varepsilon\}, \tag{5.2}$$

conflates the last $\varepsilon$ of (5.1) where it operates as a dimensionless probability with the $\varepsilon$ of (5.2) where it represents a rounding error. We refer to Davies (1993) for a discussion of this point. Other simpler metrics are also capable of dealing with rounding errors. The Kolmogorov metric is defined by

$$d_{\mathrm{ko}}(P, Q) = \sup\{|P(I) - Q(I)| : I = (-\infty, x], x \in \mathbb{R}\}. \tag{5.3}$$

Let $P_n$ be the empirical distribution of some data and $P_n^*$ be the empirical distribution of the same data after rounding. If the rounding $\delta$ is less than the minimum gap between the unrounded observations, then $d_{\mathrm{ko}}(P_n, P_n^*) = 1/n$ assuming at least one observation to have been altered. It is sometimes argued that $d_{\mathrm{ko}}$ is too weak in the data analytical sense for comparing distributions. There are stronger versions which go under the name of Kuiper metrics. The Kuiper metric of order $k$ is defined by

$$d_{\mathrm{ku},k}(P, Q) = \sup\left\{\sum_{j=1}^{k} |P(I_j) - Q(I_j)| : I_1, \ldots, I_k \text{ disjoint finite intervals}\right\}.$$
(5.4)

Kuiper metrics of order $k = 19$ are used in Davies and Kovac (2004) in the context of providing approximate densities for data. The Kolmogorov and Kuiper metrics are restricted to $\mathbb{R}$, but in higher dimensions metrics on Vapnik–Cervonenkis classes of sets retain many of their properties [see Pollard (1984)]. We refer to Davies (1993) for their use in the regression setting. The conflation of measurement error and probability in (5.1) can be avoided as follows. We define

$$(5.5) \quad d_{\mathrm{pk}}(P, Q) = \inf\{\varepsilon > 0 : P(I) \leq Q(I^\varepsilon) + \varepsilon, \text{ for all intervals } I\},$$

where $I^\varepsilon$ denotes the interval with the same center as $I$ but with length $|I|\exp(\varepsilon)$. All occurrences of $\varepsilon$ in (5.5) are now dimensionless. The idea is not new. We refer to Davies (1992, 1993).

Hampel's second argument for the Prohorov metric is that it metricizes weak convergence but we fail to see the relevance of this. The Kolmogorov metric (5.3) does not metricize weak convergence but nevertheless does have advantages over the Prohorov metric for proving central limit theorems. In particular we have

$$(5.6) \qquad d_{\mathrm{ko}}(P_n, P) = O_P(1/\sqrt{n})$$

uniformly in $P$. If $T$ is a functional with a bounded influence function $I(x, T, P)$, then under appropriate regularity conditions

$$(5.7) \quad T(P_n) - T(P) = \int I(x, T, P)\, d(P_n(x) - P(x)) + o_P(d_{\mathrm{ko}}(P_n, P)),$$

which in the light of (5.6) gives us a central limit theorem for $\sqrt{n}(T(P_n) - T(P))$. The same reasoning fails for the Prohorov metric because (5.6) does not hold [see Kersting (1978)].

**6. Metrics on $\Theta$.** We turn to the metric $D$ on $\Theta$ which quantifies the "ruinous downfall." For location in $\mathbb{R}$ the choice $D(\theta_1, \theta_2) = |\theta_1 - \theta_2|$ seems natural but the choice $|\log(\theta_1/\theta_2)|$ for scale is not quite as obvious. It does, however, have a strong justification in that numbers often have to be standardized by division by scale. If so, a scale of zero is a "ruinous downfall." In higher dimensions breakdown in scale includes the data being concentrated on a lower-dimensional hyperplane, making it impossible to identify the influence of individual covariables. Again the word breakdown would seem appropriate. In an earlier version of our paper we considered the possibility of measuring differences in the parameter $\theta$ by differences in the corresponding distributions $P_\theta$ as in $D(\theta_1, \theta_2) = d(P_{\theta_1}, P_{\theta_2})$ for an appropriate metric $d$ on the space of distributions, but this needs to be given more thought.

Tyler has pointed out that if the parameter space is compact, then the metric is bounded so that condition (3.1) of our paper cannot possibly be satisfied. This is true, but just as metrics on $\mathcal{P}$ are chosen for the problem, so we can choose metrics on $\Theta$ according to the problem. If breakdown is defined in terms of convergence to some parameter values such as those on the boundary, then we can choose an appropriate metric as follows. We start by considering the problem of scale in $\mathbb{R}$. The proof works by showing that if $\varepsilon > (1 - \Delta(P))/2$, then there exists an affine transformation $\mathcal{A}(x) = ax + b$ with $|a| \neq 1$ and, for any $n$, distributions $Q_{1n}$ and $Q_{2n}$ satisfying

$$d(P, Q_{1n}) < \varepsilon, \qquad d(P, Q_{2n}) < \varepsilon, \qquad T(Q_{1n}) = |a|^n T(Q_{2n}).$$

From this it follows that either

$$\liminf_{n \to \infty}(\min(T(Q_{1n}), T(Q_{2n}))) = 0 \quad \text{or} \quad \limsup_{n \to \infty}(\max(T(Q_{1n}), T(Q_{2n}))) = \infty.$$

Using this fact we can define the breakdown point by

$$(6.8) \quad \varepsilon^*(T, P, d, \{0, \infty\}) = \inf\{\varepsilon > 0 : \inf[T(Q) : d(Q, P) < \varepsilon] = 0$$
$$\text{or } \sup[T(Q) : d(Q, P) < \varepsilon] = \infty\}.$$

This definition makes no reference to a metric but two points on the boundary of the parameter space, 0 and $\infty$, play a special role. The metric we use in this case is $D(\theta_1, \theta_2) = |\log(\theta_1/\theta_2)|$ and, not surprisingly, the points 0 and $\infty$ also play a special role here. The result is that $\varepsilon^*(T, P, d, \{0, \infty\}) = \varepsilon^*(T, P, d, D)$. We see that breakdown as defined by (6.8) can be reformulated in terms of an appropriately chosen metric on the parameter space $\Theta$. This remains true even if $\Theta$ is compact. Suppose $\Theta$ is equipped with a metric $D^*$, bounded or not, and that some parameter value $\theta_0$ is regarded as breakdown, for example, 0 in the scale context or 1 in the correlation context. We define the metric $D$ on $\Theta$ by

$$(6.9) \qquad D_{\theta_0}(\theta_1, \theta_2) = \left| \frac{1}{D^*(\theta_1, \theta_0)} - \frac{1}{D^*(\theta_2, \theta_0)} \right|.$$

It follows that if we keep $\theta_1$ constant, then $D_{\theta_0}(\theta_1, \theta_2)$ tends to infinity if and only if $\theta_2$ tends to $\theta_0$. If we define in analogy to (6.8)

$$(6.10) \quad \varepsilon^*(T, P, d, \{\theta_0\}) = \inf\{\varepsilon > 0 : \inf[D^*(\theta_0, T(Q)) : d(Q, P) < \varepsilon] = 0\},$$

then clearly $\varepsilon^*(T, P, d, \{\theta_0\}) = \varepsilon^*(T, P, d, D_{\theta_0})$. If there is a set of parameter values $\Theta_0$ which are regarded as breakdown, for example, the boundary points, we define

$$(6.11) \qquad D_{\Theta_0}(\theta_1, \theta_2) = \sup\{D_{\theta_0}(\theta_1, \theta_2) : \theta_0 \in \Theta_0\}$$

and again we have a metric which can be used to define breakdown. We define

$$(6.12) \quad \varepsilon^*(T, P, d, \Theta_0)$$
$$= \inf\{\varepsilon > 0 : \inf\{\inf[D^*(\theta_0, T(Q)) : d(Q, P) < \varepsilon] : \theta_0 \in \Theta_0\} = 0\}$$

and it follows that $\varepsilon^*(T, P, d, \Theta_0) = \varepsilon^*(T, P, d, D_{\Theta_0})$ and also

(6.13) $\qquad \varepsilon^*(T, P, d, \Theta_0) = \inf\{\varepsilon^*(T, P, d, \{\theta_0\}): \theta_0 \in \Theta_0\}.$

Grize (1978), as we shall see below, defines breakdown as the minimum contamination such that all points in $\Theta_0$ are reachable and not just some such point. This can be accommodated by defining

(6.14)
$$\varepsilon^{**}(T, P, d, \Theta_0) \\ = \inf\{\varepsilon > 0 : \inf\{\theta \in \Theta_0 : \sup[D_{\theta_0}(T(P), T(Q)) : d(P, Q) < \varepsilon] = \infty\}\}.$$

In contrast to (6.13) this definition results in

(6.15) $\qquad \varepsilon^{**}(T, P, d, \Theta_0) = \sup\{\varepsilon^*(T, P, d, \{\theta_0\}): \theta_0 \in \Theta_0\}.$

There are no doubt other variations. The conclusion is that if breakdown is defined as convergence to some set of exceptional parameter values, then this can be described by a metric as required in our theorem. It still leaves open the question as to whether such a definition of breakdown is sensible but this can only be answered on a case-to-case basis.

**7. Breakdown point?** Rousseeuw in his contribution argues for the use of breakdown "value" rather than "point." We do not quite understand his reasoning and while usage is never absolute, we do not see any advantage in replacing "point" by "value." Hampel mentions the analysis of variance as one situation where the term breakdown point may not be appropriate. In the simple two-way table breakdown occurs if the majority of observations in any one row or column are badly contaminated, but this is too pessimistic and gives an artificially low breakdown point. In Terbeck and Davies (1998) the breakdown or interaction *patterns* for the two-way table are characterized and it is shown how these are related to the $L_1$-solution and to Tukey's median polish. Other articles concerned with patterns are Ellis and Morgenthaler (1992) for $L_1$ regression and Kuhnt (2000) for contingency tables.

**8. Affine equivariance.** We agree with Hampel that affine equivariance is not always a requirement and in two or more dimensions it is more difficult to justify than in one. We made comments to this effect in our paper. Nevertheless it is not always the case that outliers are apparent in the single coordinates and to find these some sort of equivariance would seem to be required. An example of a simple data set for which it is not sufficient just to look at the coordinates is given on page 57 of Rousseeuw and Leroy (1987). Programs based on high breakdown methods are now readily available and in our opinion should be used in a routine manner [Becker and Gather (1999, 2001), Rocke (1996) and Rousseeuw and van Driessen (1999)]. The costs are negligible and the returns can be substantial.

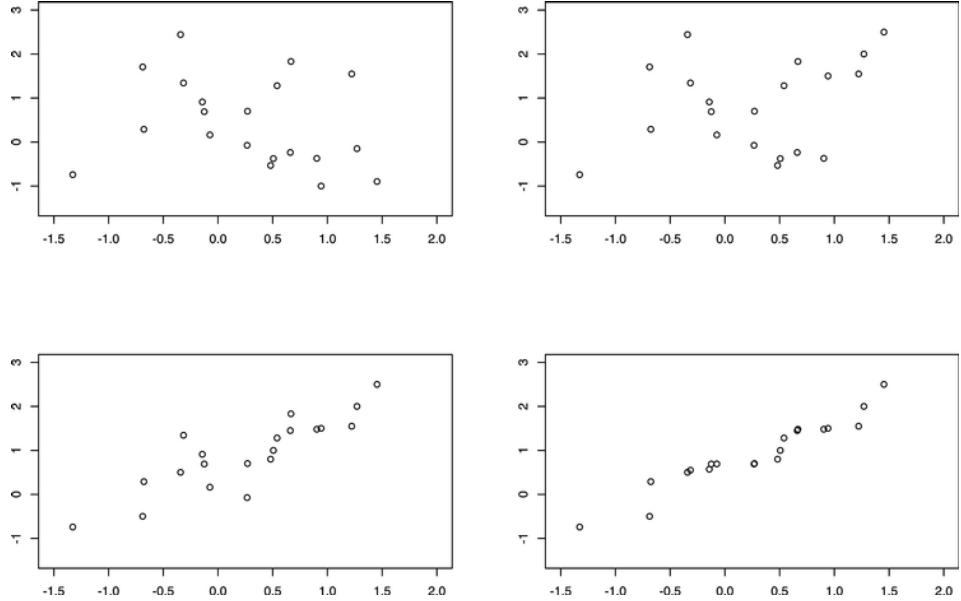

FIG. 1. *Samples differing from the initial sample* (upper left) *by* 3 *points* (upper right), 9 *points* (lower left) *and* 14 *points* (lower right) *with the rank-correlation changing from* −0.332 *to* 0.278, 0.878 *and* 1, *respectively.*

**9. Correlation.** This brings us to the perhaps most important part of the discussion. Hampel argues strongly that correlation provides an example of a useful concept of breakdown which does not have an equivariance structure. We argue that he is wrong on both counts: the concept is not useful and it does have an equivariance structure, albeit a simple one. We give a detailed reply which touches on many of the points discussed so far. Grize (1978) gives two definitions of breakdown for a rank correlation functional $T_{\rm rc}$. The first reads [see (6.14), (6.15)]

$$(9.16) \quad \varepsilon^{**}(T_{\rm rc}, P, d, \{-1, 1\}) = \inf\{\varepsilon > 0 : \sup\{T_{\rm rc}(Q) : d(P, Q) < \varepsilon\} = 1,$$
$$\inf\{T_{\rm rc}(Q) : d(P, Q) < \varepsilon\} = -1\}$$

and the second reads

$$(9.17) \quad \varepsilon^{**}(|T_{\rm rc}|, P, d, \{0, 1\}) = \inf\{\varepsilon > 0 : \sup\{|T_{\rm rc}(Q)| : d(P, Q) < \varepsilon\} = 1,$$
$$\inf\{|T_{\rm rc}(Q)| : d(P, Q) < \varepsilon\} = 0\}$$

for some appropriate metric $d$. For the total variation metric Grize calculates the breakdown points of Kendall's and Spearman's rank correlation for some particular distributions. We carry out a small experiment for Spearman's functional $T_{\rm sc}$. The top left panel of Figure 1 shows 20 data points with $T_{\rm sc}(P_n) = -0.332$. Initially there are various sets of six points for which $y_i = h(x_i)$ with $h$ a nondecreasing function. We choose one and then move one of the remaining points at a time until finally after 14 moves all the points satisfy $y_i = h(x_i)$ with $h$ nondecreasing. For the final sample the rank correlation is 1 and we have, according to (9.16), breakdown. The top right panel of Figure 1 shows the sample after three moves, the bottom left after nine moves. The final sample is shown in the bottom right panel. At no stage do we experience a breakdown. Each sample differs only slightly from the previous one and the values of the rank correlation are perfectly reasonable for the sample they refer to. In Hampel's terminology there is no pole. In a similar vein, Figure 2 shows a distribution considered by Grize for which he calculates the breakdown point 0.1 of Spearman's rank correlation in the sense of (9.17). The top left panel shows the initial distribution and the top right panel the same data after a monotone transformation. A breakdown to zero is shown in the bottom left panel and to 1 in the bottom right panel. In our opinion the bottom right panel is the only one where one would not a priori question any observation and yet this is classified as breakdown. We now play a similar game in one dimension and consider a simple standard normal sample of size 20. We consider the median and as breakdown corresponds to an arbitrarily large value of the median we start with 100. The game is now to alter the initial sample point by point until after ten moves the value of the median is at least 100. The moves are almost prescribed. We choose any observation from the original sample and move it about 200 units to the right. After ten moves the median assumes a value of about 100. There is no other strategy. Even the first move alters the sample in a manner which distinguishes it immediately from the initial sample. Furthermore, when we progress from the ninth to the tenth move the median suddenly jumps from a value of about zero to one of about 100. We think this situation can be described by the word "breakdown." Moreover, it holds for any translation equivariant functional if one replaces the points carefully as in (6.2) and not as in (6.3). Any such functional must break down by the tenth move at the latest.

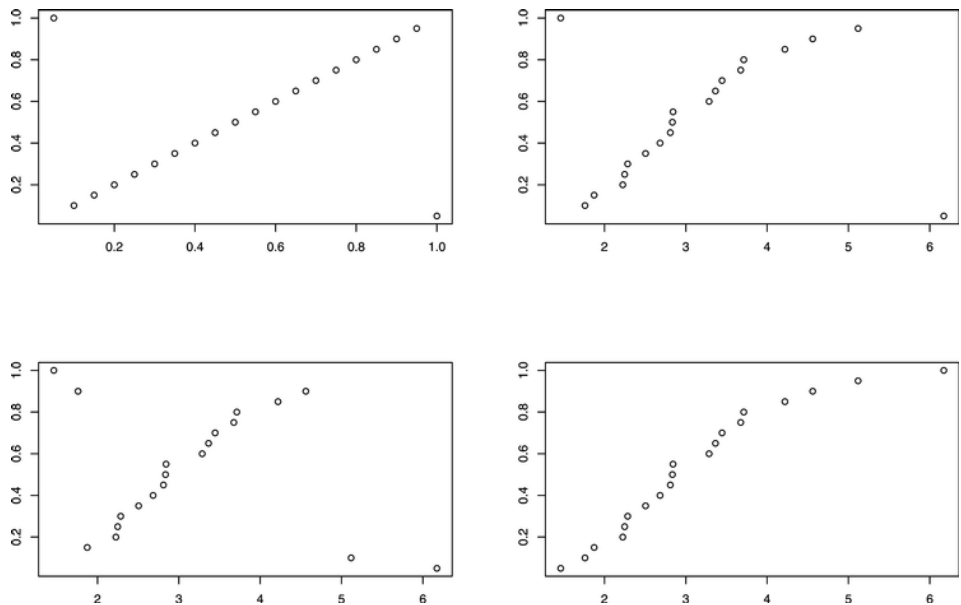

FIG. 2. *The upper left panel shows a distribution considered by Grize (1978). The upper right panel shows the same data after a monotone transformation. The bottom left panel shows the breakdown [in the sense of (9.17)] of Spearman's rank correlation to zero. The bottom right panel shows the breakdown of Spearman's rank correlation to 1.*

We now consider the usual linear correlation functional $T_{lc}$. For the initial data set of Figure 1 its value is $-0.258$. If we take any observation and move it to the point $(\gamma, \gamma)$ and let $\gamma$ tend to infinity, then $T_{lc}$ tends to 1. In this situation it seems reasonable to use the word breakdown but perhaps discontinuity would be a better description. We analyze the problem more closely. Linear correlation can be placed into the context of our paper by introducing the following group structure. We define $G_{lc}$ to be the group of transformations $g: \mathbb{R}^2 \to \mathbb{R}^2$ with

$$(9.18) \qquad g(x, y) = (a_1 x + b_1, a_2 y + b_2)$$

with $a_1 a_2 \neq 0$. An equivariant functional $T$ is one which satisfies

$$(9.19) \qquad T(P^g) = \operatorname{sgn}(a_1 a_2) T(P).$$

Clearly the usual linear correlation functional $T_{lc}$ is equivariant w.r.t. this group. The metric on the space of distributions is taken to be the strip metric

$$(9.20) \qquad d_{\mathcal{ST}}(P, Q) = \sup\{|P(C) - Q(C)| : C \in \mathcal{ST}\},$$

where $\mathcal{ST}$ denotes the set of strips $C$

$$(9.21) \qquad C = \{(x, y) : -\delta \leq ax + by + c \leq \delta; a, b, c \in \mathbb{R}, \delta \in \mathbb{R}_+\}.$$

We note that this metric is also "correct" as it is invariant under the group $G_{lc}$:

$$(9.22) \qquad d_{\mathcal{ST}}(P, Q) = d_{\mathcal{ST}}(P^g, Q^g), \qquad g \in G_{lc}.$$

There is also a version of this metric which corresponds to (5.5) [see Davies (1993)]. To fit into the structure in our paper we also require a metric $D$ on the parameter space $\Theta = [-1, 1]$. The precise metric is not important because of the simple nature of the equivariance given in (9.19). To be concrete we put in (6.9)

$$D^*(\theta_1, \theta_2) := |\tan(\theta_1 \pi/2) - \tan(\theta_2 \pi/2)|,$$

which is consistent with the desire to have breakdown at $\pm 1$.

From (9.19) we see that the condition $G_1 \neq \varnothing$ is not satisfied and Theorem 3.1 does not provide a nontrivial upper bound. Indeed, there is an equivariant correlation functional with breakdown point 1, namely $T_{lc}^o \equiv 0$, but to forestall protests we give another. For an empirical distribution $P_n$ we define

$$(9.23) \qquad T_{lc}^*(P_n) = \frac{1}{N} \sum_{I, |I| \geq 3} T_{lc}(I),$$

where $I$ is a subset of the data containing $|I|$ observations, $N = 2^n - n - n(n-1)/2$ and $T_{lc}(I)$ is, by an abuse of notation, $T_{lc}$ evaluated at the empirical measure based on the set of observations in $I$. The functional $T_{lc}^*$ is equivariant and also Fisher consistent. To calculate the breakdown point we consider an empirical measure $P_n$ deriving from a sample of size $n$ from a continuous distribution on $\mathbb{R}^2$ and another empirical measure $Q_n$. We assume that the supports of each are contained in some compact set $K$. The reason for these assumptions is to reduce complications due to the fact that the linear correlation coefficient as usually defined requires the existence of moments. We consider a sequence of $Q_n$ with $\lim_{n \to \infty} T_{lc}^*(Q_n) = 1$. From (9.23) it follows that the support of $Q_n$ must be contained in a strip

$$C_n = \{(x,y) : -\delta_n \leq a_n x + b_n y + c_n \leq \delta_n\}$$

with $\lim_{n \to \infty} \delta_n = 0$. As $P_n(C_n) \leq 2/n$ for sufficiently large $n$ we have $d_{\mathcal{ST}}(P_n, Q_n) \geq 1 - 2/n$ and hence

$$(9.24) \qquad \varepsilon^{**}(T_{lc}^*, P_n, d_{\mathcal{ST}}, \{-1, 1\}) \geq 1 - 2/n$$

for this class of probability measures. We generalize this result in a manner which emulates the setting of our paper. As $G_1$ is empty we reformulate the definition of the functional $\Delta(P)$ of (3.3) in our paper as follows. We set

$$(9.25) \quad \Delta(P) = \sup\{P(B) : T(Q) \text{ not definable for } Q \text{ with } \mathrm{supp}(Q) \subset B\}.$$

For example, in the case of scale in $\mathbb{R}$ the relevant sets $B$ are singletons and a measure concentrated on a singleton must have scale either zero or $\infty$ to be equivariant, both of which are excluded. If, following Grize, a linear correlation of $\pm 1$ is defined to be breakdown, the corresponding sets are lines and this leads to

$$(9.26) \quad \Delta_+(P) = \sup\{P(C) : C = \{(x,y) : ax + by + c = 0\}, ab \leq 0\},$$

$$(9.27) \quad \Delta_-(P) = \sup\{P(C) : C = \{(x,y) : ax + by + c = 0\}, ab \geq 0\},$$

and it follows that

$$(9.28) \qquad \varepsilon^{**}(T_{lc}^*, P, d_{\mathcal{TU}}, \{-1, 1\}) = 1 - \min\{\Delta_+(P), \Delta_-(P)\}.$$

The reasoning can be extended to rank correlation and this gives a more elegant theory as there are no problems with moments. The appropriate group is $G_{\rm rc}$ which consists of all transformations $g:\mathbb{R}^2 \to \mathbb{R}^2$ of the form

$$(9.29) \qquad g((x,y)) = (\zeta(x), \eta(y)), \qquad \zeta, \eta, \mathbb{R} \to \mathbb{R},$$

where each of $\zeta$ and $\eta$ is either strictly increasing or strictly decreasing. A correlation functional $T_{\rm rc}$ is equivariant with respect to this group if

$$(9.30) \qquad T_{\rm rc}(P^g) = \text{sgn}(\zeta \circ \eta) T_{\rm rc}(P),$$

where $\text{sgn}(\zeta) = \pm 1$ depending on whether $\zeta$ is strictly increasing or decreasing. The natural metric is the tube metric

$$(9.31) \qquad d_{\mathcal{TU}}(P,Q) = \sup\{|P(C) - Q(C)| : C \in \mathcal{TU}\},$$

where $\mathcal{TU}$ denotes the set of monotonic tubes $C$

$$(9.32) \quad C = \{(x,y) : -\delta \leq h(x) + y \leq \delta, h : \mathbb{R} \to \mathbb{R} \text{ strictly monotonic}, \delta \in \mathbb{R}_+\}.$$

The metric is "correct" in that it is invariant with respect to the group $G_{\rm rc}$:

$$(9.33) \qquad d_{\mathcal{TU}}(P,Q) = d_{\mathcal{TU}}(P^g, Q^g), \qquad g \in G_{\rm rc}.$$

As we now require correlations of $\pm 1$ only for data points which are strictly increasing or decreasing, we define analogously to (9.26) and (9.27),

$$(9.34) \quad \Delta_+(P) = \sup\{P(C) : C = \{(x,y) : y = h(x)\}, h \text{ strictly increasing}\},$$

$$(9.35) \quad \Delta_-(P) = \sup\{P(C) : C = \{(x,y) : y = h(x)\}, h \text{ strictly decreasing}\}.$$

From this it follows for Spearman's rank correlation functional $T_{\rm sc}$ that

$$(9.36) \qquad \varepsilon^{**}(T_{\rm sc}, P, d_{\mathcal{TU}}, \{-1, 1\}) = 1 - \min\{\Delta_+(P), \Delta_-(P)\}.$$

In fact (9.36) holds for any functional for which $T_{\rm rc}(P) = 1$ or $-1$ if and only if $\Delta_+(P) = 1$ or $\Delta_-(P) = 1$, respectively, and consequently it also holds for Kendall's $\tau$. The appearance of min in (9.28) and (9.36) is due to Grize's definition (9.16) of breakdown which refers to both boundary points. Usable estimates of $\Delta_+(P_n)$ are available for empirical measures $P_n$ deriving from nonatomic i.i.d. random variables in each component; that is, the components are also independent. Let the sample be $(X_i, Y_i)$, $i = 1, \ldots, n$, and consider the points $(X_{i1}, Y_{i1}), \ldots, (X_{ik}, Y_{ik})$ with the $X_{ij}$, $j = 1, \ldots, k$, in increasing order. The points lie on some curve $y = h(x)$ for a strictly increasing $h$ if and only if the $Y_{ij}$, $j = 1, \ldots, k$, are also in increasing order. The probability of this is $1/k!$. There are $\binom{n}{k}$ different samples of size $k$ and we see that the probability that at least $k$ points lie on some curve $y = h(x)$ is at most

$$\frac{1}{k!}\binom{n}{k} \leq \frac{n^k}{(k!)^2}.$$

By maximizing over $k$ we obtain

$$\Delta_+(P_n) = O(1/\sqrt{n}\,)$$

and it follows from (9.36)

(9.37) $$\varepsilon^{**}(T_{\text{sc}}, P_n, d_{\mathcal{TU}}, \{-1, 1\}) \geq 1 - O(1/\sqrt{n}\,).$$

The fact that (9.36) also holds for Kendall's $\tau$ apparently contradicts Hampel's comments, but this is not so because it is definition (9.17) of breakdown to which Hampel's comments apply. To proceed we consider the problem of maximizing $\Delta_+(P)$ subject to $T_{\text{rc}}(P) = 0$. For $T_{\text{sc}}$ the answer is $\Delta_+(P) = \sqrt[3]{1/2}$, which is attained at a distribution for which the rank of $x_i$ is $i$ and the rank of $y_i$ is $k+i$, $1 \leq i \leq n-k$, and $n-i+1$, $k+1 \leq i \leq n$, with $k = n\sqrt[3]{1/2}$. The corresponding result for Kendall's $\tau$ replaces $\sqrt[3]{1/2}$ by $\sqrt{1/2}$. If now $Q$ is any distribution with $\Delta_+(Q) = 1$, it follows that the breakdown point [in the sense of (9.17)] at $Q$ is $1 - \sqrt[3]{1/2} = 0.2063$ for Spearman's rank correlation and $1 - \sqrt{1/2} = 0.2929$ for Kendall's $\tau$. If we now move only half the mass of $1 - \sqrt[3]{1/2}$, it is clear that we can obtain distributions $Q_1$ and $Q_2$ at which Spearman's and Kendall's rank correlations have breakdown points of $(1 - \sqrt[3]{1/2})/2$ and $(1 - \sqrt{1/2})/2$, respectively, and that these are the smallest possible breakdown points. We have not understood Hampel's claim $BP(K) = \frac{3}{2}BP(S)$ as, as far as we can see, these refer to different distributions, one with $\Delta_+(Q_1) = 0.85$ and one with $\Delta_+(Q_2) = 0.9$, but this is only a minor point. The tube metric $d_{\mathcal{TU}}$ is stronger than the strip metric $d_{\mathcal{ST}}$ but considerably weaker than the total variation metric $d_{\text{tv}}$ used by Grize. In particular it allows for wobbling of the observations. We also note that neither metric suffers from the deficiency of the Prohorov metric of mixing dimensionless probabilities with measurement units. The class $\mathcal{ST}$ of strips has polynomial discrimination but not the class $\mathcal{TU}$ as there are arbitrarily large finite subsets of $\mathbb{R}^2$ which can be shattered by $\mathcal{TU}$ [see Pollard (1984)].

Finally, we point out the differences to Theorem 3.1. If we take the usual definition of breakdown as a worst-case situation rather than Grize's definition which is a sort of best worst-case definition, then the breakdown point of Kendall's or Spearman's rank correlation is

(9.38) $$\varepsilon^*(T_{\text{lc}}^*, P, d_{\mathcal{TU}}, \{-1, 1\}) = 1 - \Delta(P),$$

where

$$\Delta(P) = \max\{\Delta_+(P), \Delta_-(P)\}.$$

In the general situation we argue as follows. Let $\mathcal{P}_0$ denote the set of distributions at which $T$ breaks down, which means that their support is contained in some exceptional subset of the sample space as in (9.25). Suppose $0 < \Delta(P) < 1$ and choose an exceptional subset $B_0$ of the sample space with $P(B_0) = \alpha$, $0 < \alpha < \Delta(P)$. If we define $Q_0(\cdot) = P(\cdot \cap B_0)/P(B_0)$ and $Q_1(\cdot) = P(\cdot \cap (\mathcal{X} \setminus B_0))/(1 - P(B_0))$, then $Q_0$ and $Q_1$ are probability measures with $P = \alpha Q_0 + (1-\alpha) Q_1$. If the metric on $\mathcal{P}$ satisfies (2.2) of our paper, we see that $d(P, Q_0) \leq 1 - \alpha$ and this implies

$$\varepsilon^*(T, P, d, D) \leq 1 - \Delta(P). \tag{9.39}$$

This differs from the claim of Theorem 3.1 by the factor of $1/2$ and it is precisely the group structure which produces this factor. Because of equivariance things start going wrong before one reaches an arbitrarily small neighborhood of some point in $\mathcal{P}_0$. As Tyler mentions in his contribution, heuristic justifications for the factor of $1/2$, such as not being able to distinguish between good and bad data, are too vague. One of the challenges of this paper is to obtain the factor of $1/2$ or even some other factor without an equivariance structure.

**10. Principal component vectors.** Tyler argues that it may be possible to define a reasonable concept of breakdown for principal component vectors, although he recognizes that there are problems involved. The idea is that breakdown occurs if contamination results in the first principal component vector being orthogonal to the first principal component vector without contamination. This example cannot be reformulated in terms of metrics as described in Section 6, as there is no special set of parameter values $\Theta_0$. Furthermore, it is not possible to adjust the proof of Theorem 3.1 to include this case. However, we shall now argue that it does not make sense to talk about the breakdown of the principal component vectors without reference to the corresponding eigenvalues.

Consider a two-dimensional data set for which the eigenvalues are the same. The set of first principal component vectors is now the set of all points on the unit circle with $\theta$ and $-\theta$ being identified. The smallest alteration of any observation will cause the space to collapse to a single direction, say $\theta_1 = (1,0)$, with the second principal component vector $\theta_2 = (0,1)$ being orthogonal to it. It is clear that there exists an arbitrarily small perturbation of the original data set such that $\theta_1 = (0,1)$ and $\theta_2 = (1,0)$. In other words, there exist data sets for which arbitrarily small perturbations cause breakdown in the first principal component vector. The perturbations can be so small as to be nondetectable and any computer of finite precision or some nonoptimal numerical recipe may result in the wrong answer and be the "cause" of the breakdown. It seems to us that this situation is one which is not describable by the word "breakdown." We cannot think of any useful statistical procedure which can be made to break down by the smallest of perturbations of the data set. In practice, of course, use is not just made of the first principal component vector but of all those principal component vectors for which the eigenvalues are in some sense large. One data analytical strategy is to look at the two-dimensional plots on the first two principal component vectors, and here it is irrelevant if they are the wrong way round. The principal component vectors are defined as those directions where the variability of the data, however measured, is particularly large. The word "breakdown" can be more appropriately applied to a situation in which the large variability is the result of outliers and causes a direction of small variability to become one of large variability. It seems to us to be clear that the size of the eigenvalues will have to be taken into account. Principal component vectors do not therefore constitute a counterexample to our meta claim of no nontrivial theory of breakdown without groups. In spite of this Tyler has alerted us to the possibility of breakdown being defined in terms of a relationship between two parameter values rather than closeness to some specific parameter values. We cannot exclude the possibility of there being some perfectly reasonable concept of breakdown of this nature.

**11. Fisher consistency.** Hampel does not like the example of regression through the origin and neither do we. It was included as an answer to a referee as to whether it was possible to construct a Fisher consistent functional with a breakdown point of at least $1/3$. Fisher consistency seems to be the obvious candidate to replace group equivariance as a desirable property of a functional. We do not give a theorem as there are difficulties in defining what is meant by a reasonable parametric family, but a parametric family $\{P_\theta, \theta \in \Theta\}$ typically forms a very sparse subset of the set of all models. This is indicated by Figure 3 where the line represents the family of models in the space of all probability measures and the circles indicate an infinitesimal neighborhood. Fisher consistency describes the behavior of the functional only in the infinitesimal neighborhood. We are left free to define the functional elsewhere and this is what we exploit in our example. Equivariance considerations prevent this form of local definition. The orbits connect points which are far apart in the space of probability models and this prevents constructions such as the one we give.

**12. The samples (6.2) and (6.3).** Hampel, Rousseeuw and Tyler all comment on the samples (6.2) and (6.3). Rousseeuw correctly remarks that one can often calculate the breakdown point of a functional directly and that such direct proofs do not rely on a repetition. Hampel says, also correctly, that the unnamed functional (there are many) must have a low breakdown and suggests that perhaps some small print is missing. What is missing is some large print explaining exactly what we intended with these two examples. Tyler saw clearly what was intended and has made some very interesting comments on (6.2) and (6.3). He also explicitly mentions the connection with the area of computer vision which was one of our motivations as we indicate below. He points out that such apparently well-understood functionals such as appropriately tuned redescending $M$-functionals can exhibit the same behavior. We confess to not having been aware of this and we would have chosen another example had we known. As Tyler says, under appropriate conditions redescending affine equivariant $M$-functionals do not break down even under 99% contamination. This is exactly the phenomenon to which we intended to bring attention.

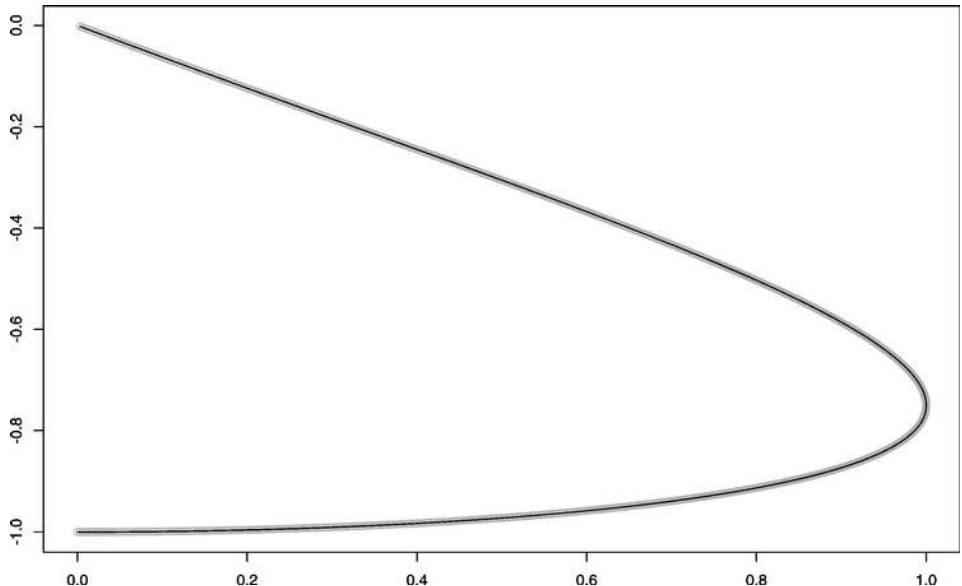

Fig. 3. *A thin parametric model and an infinitesimal neighborhood within which Fisher consistency becomes relevant.*

The proof of Theorem 3.1 relies in part on exactly reproducing a portion of the data elsewhere. If there is no exact repetition as in (6.3), there will be many equivariant functionals which do not break down. One choice for sample (6.3) is

$$(12.40) \quad (T_l(\mathbf{x}_n), T_s(\mathbf{x}_n)) = \arg\min_{\mu,\sigma} \left\{ \sum_{j=1}^{3} r_{(i)}(\mu,\sigma)^2 \right\},$$

where

$$r_i(\mu,\sigma)^2 = \min_{1 \leq j \leq n} \left( \frac{x_j - \mu}{\sigma} - z_i \right)^2$$

and $z_1 = 1.5, z_2 = 1.8$ and $z_3 = 1.3$. In this connection we mention Oja's example of linear regression at the end of his Section 1. We fail to follow his argument as to why the estimate becomes uninformative. As it stands, the argument seems to make no use of the assumption $n = 2k$, in which case we can put $k = 1$ and the conclusion would seem to be that every regression equivariant functional has a breakdown point of $1/n$. If $n = 2k$ is implicitly meant, then breakdown occurs only if we cannot distinguish between the two samples. If we can distinguish between the two samples, for his example if $x_1 = x_2 = \cdots = x_k = 0$, then what is claimed as breakdown is nothing but equivariance (see Section 2 above).

At first glance the functional (12.40) may seem very artificial but this is not so. It is constructed to find a particular pattern in the sample, namely affine transformations of 1.5, 1.8 and 1.3. Figure 4 shows the smile of the Cheshire cat and the problem is to locate it in a sea of noise into which it gradually disappears. This is only possible as the noise does not reproduce the signal. For real examples from the area of computer vision we refer to Wang and Suter (2004). In analytical terms one can define a modified breakdown point by

$$(12.41)\ \varepsilon^*(T, P, d, \mathcal{H}) = \inf\{\varepsilon > 0 : \sup |T(P) - T(Q)|, d(P, Q) < \varepsilon, Q \in \mathcal{H}\},$$

where $\mathcal{H}$ specifies what you want to protect yourself against. If $\mathcal{H}$ does not allow a repetition of the signal elsewhere, then affine equivariant functionals can attain breakdown points higher than $1/2$. Moreover, the usual high breakdown functionals are typically of no help in this situation [see Wang and Suter (2004)].

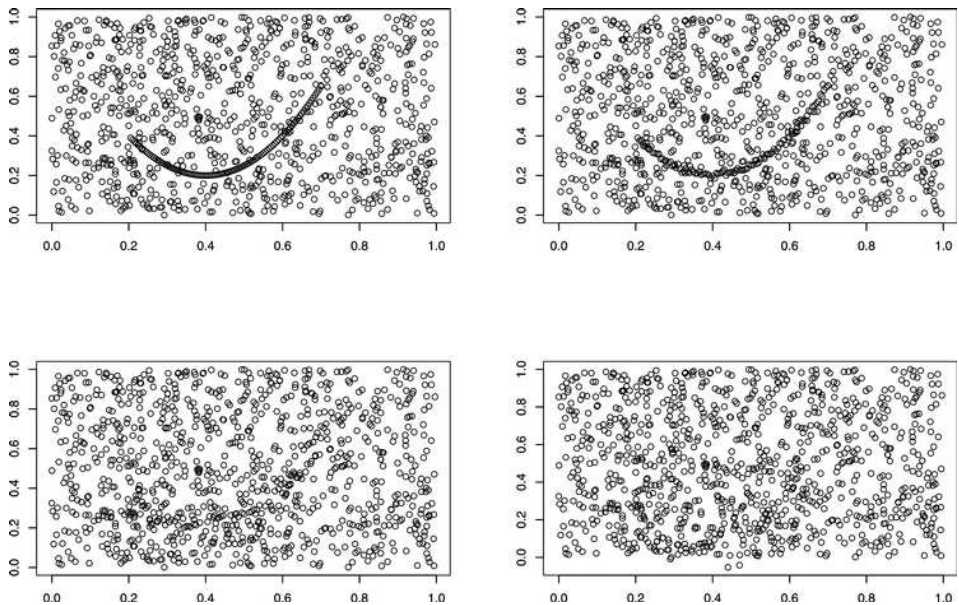

Fig. 4. *The smile of the Cheshire cat gradually disappearing as it is corrupted by noise.*

**13. Nonparametric statistics.** In this paper we have shown that the concept of breakdown point has been generally accepted only in situations where there is a group structure sufficiently rich to allow the calculation of a nontrivial upper bound for the breakdown point. In his contribution Hampel speculates that this could be the reason why the breakdown point for correlation coefficients has not yet been widely accepted. In spite of this and as pointed out by an Associate Editor, we have not proved a theorem to the effect that a breakdown point is only sensible when a rich group structure exists. It is difficult to imagine what such a theorem would look like. Nevertheless the paper, the contributions of the discussants and our reply do seem to indicate that it will not be easy to come to an acceptable definition of breakdown with a nontrivial upper bound without a group structure. There is perhaps another reason why some definitions of breakdown have been successful. They are defined for so-called nonparametric functionals in the sense of Bickel and Lehmann (1975a, b). One can always calculate the median of a distribution in $\mathbb{R}$ and this is not associated with a restrictive stochastic model. We wish to emphasize this as we have the impression that it is sometimes assumed that functionals are only to be applied to data which is generated by some stochastic model but with contamination. Genton and Lucas entitle a section "Breakdown point for (in)dependent observations," which suggests at least to us that they distinguish between samples which are generated by independent random variables and those which are not. The title of Genton and Lucas (2003) also tends in this direction. They write "$\mathcal{Y}$ is the set of all allowable samples" and later "$\mathcal{Y}$ is the set of all stationary AR(1) processes." We think the intention is clear. The data are generated by a stationary AR(1) process and then contaminated by the outliers. On the other hand, the only possible mathematical interpretation of "the set of all stationary AR(1) processes" is the support of the model. As the support of an AR(1) process with Gaussian innovations is $\mathbb{R}^n$, this means simply all samples, $\mathcal{Y} = \mathbb{R}^n$. Thus what at first glance seems plausible turns out to be untenable. This is the reason why the restrictions we place on data are analytical ones and not distributional ones. The median can be successfully applied to data which are very obviously dependent [see, e.g., Davies, Fried and Gather (2004)], but consider the data in Figure 5. For these data it makes no sense to artificially restrict $\rho$ to the interval $[-1, 1]$. Rather one would point out that the data can be well approximated by an AR(1) model with $\rho = -1.25$ but not by a stationary AR(1) model. This simple message is lost if we are forced to specify a value in $[-1, 1]$. If breakdown is going to be meaningful in such a situation we suspect that it should be applied to a statistical procedure and not to the behavior of a single functional.

**14. Breakdown without groups and alternatives.** We have argued above that the only generally accepted definitions of breakdown are in situations where there is a sufficiently rich group and equivariance structure. If a need is felt to extend it to other situations, we state what we think are the minimal requirements. First, the definition should be capable of being made precise. He argues that breakdown is the smallest fraction of contamination which makes a test statistic "uninformative or unusable." Later he argues that breakdown should have the same degree of vagueness as he claims to be the case with outliers. He continues that "when every statistician starts to talk about his or her own notion of a breakdown point, I think we have made it." We think there are dangers in such an attitude. Ostensive definitions of breakdown with statisticians pointing in all directions are unlikely to contribute to a general acceptance of the word. Intuition is important, but just as is the case with outliers [see Davies and Gather (1993)] much is to be gained by undertaking the attempt to give a precise definition and to investigate its consequences. This not only deepens the understanding, it also sharpens the intuition. Semantics is important and we think that any generalizations of the concept of breakdown should be such as to be recognizably referring to some common element, in particular the presence of some natural pole. Second, it is not sufficient to give a definition of breakdown and show that it gives the correct answer in some particular cases. A definition of breakdown should be subjected to some form of analysis, including its stability under perturbations. The onus is on those who propose definitions of breakdown to do this. Third, the definition should be simple and intuitively appealing. Here we agree completely with He. If it requires more than sixty seconds to understand a definition, it is probably bad. Fourth, when calculating breakdown points use should be made of metrics which can accommodate rounding errors. Gross error neighborhoods and the total variation metric are too restrictive. Fifth, the definition should not be too restrictive and only apply to one single functional. It should apply to a whole family of functionals which offer different possibilities of quantifying the feature of interest, location, scale, correlation or nonparametric function. Sixth, there should be a class of reasonable functionals for which it makes sense to compare breakdown points. If such a definition of breakdown is not possible, there are alternatives. One is simply to compare different functionals by their continuity or bias properties, again if possible in weak metrics. For this to make sense it is not necessary that an explosion occurs. It may be that this proves more useful than trying to extend the idea of breakdown to situations for which it is not suitable.

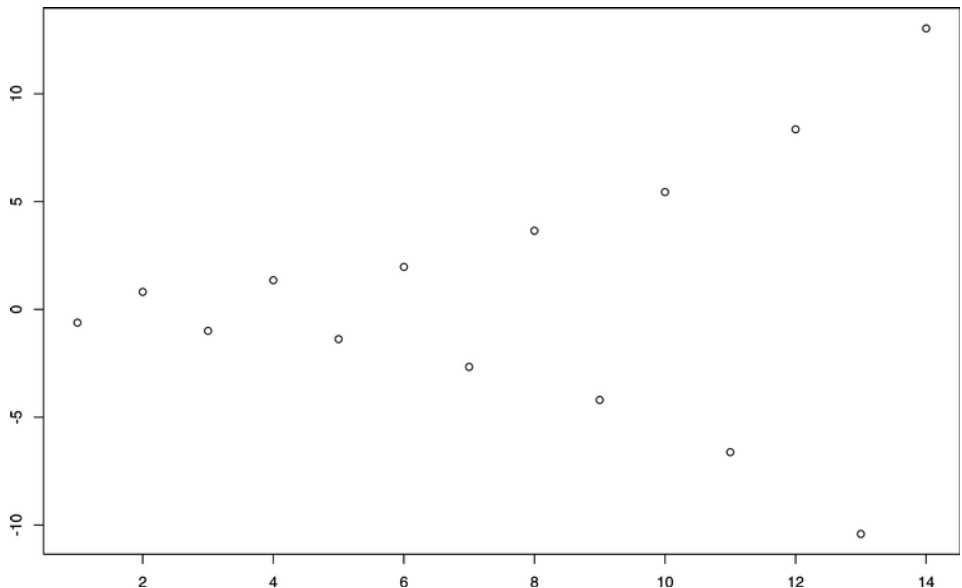

FIG. 5. *A sample of size* 14 *generated by* $X_{t+1} = -1.25X_t + 0.2Z(t)$ *with* $Z(t)$ *standard Gaussian white noise.*

**15. Conclusion.** We thank all discussants for their contributions and hope that the disagreements that are apparent have been clarified by our rejoinder. In our paper we have not proved that breakdown without equivariance is not a sensible concept. On the other hand, in all situations we are aware of in which there is no or little equivariance (made precise by our main theorem), then either (i) breakdown points of 1 are attainable or (ii) the word breakdown is inappropriate (the movement from the top right to the bottom right panel of Figure 2) or (iii) the very definition of breakdown point is inadequate. An example without these weaknesses would be interesting.

Fachbereich 06—Mathematik
und Informatik
Universität Duisburg–Essen
45117 Essen
Germany
e-mail: davies@stat-math.uni-essen.de

Fachbereich Statistik
Universität Dortmund
44221 Dortmund
Germany
e-mail: gather@statistik.uni-dortmund.de